\documentclass[10pt,leqno]{amsart}
\usepackage{amssymb}
\usepackage{amscd}
\usepackage{amsmath}
\usepackage{amsthm}

\newcommand{\Q}{\mathbb Q}
\newcommand{\Z}{\mathbb Z}
\newcommand{\R}{\mathbb R}
\newcommand{\C}{\mathbb C}
\newcommand{\G}{\mathcal G}
\newcommand{\F}{\mathbb F}

\newcommand{\A}{\mathbb A}

\makeatletter
\def\l@section{\@tocline{1}{4pt}{1pc}{}{}}
\def\l@subsection{\@tocline{2}{0pt}{2pc}{5pc}{}}
\makeatother

\begin{document}
\title{Siegel modular forms of genus $2$ attached to elliptic curves}
\author{Dinakar Ramakrishnan and Freydoon Shahidi}
\date{}

\maketitle

\medskip

\section*{Introduction}

\bigskip

The object of this article is to construct certain classes of {\it
arithmetically significant, holomorphic Siegel cusp forms $F$ of
genus $2$}, which are {\it neither of Saito-Kurokawa type}, in which
case the degree $4$ spinor $L$-function $L(s, F)$ is divisible by an
abelian $L$-function, {\it nor of Yoshida type}, in which case
$L(s,F)$ is a product of $L$-series of a pair of elliptic cusp
forms. In other words, for the holomorphic cusp forms $F$ we
construct below, which includes those of scalar weight $\geq 3$, the
cuspidal automorphic representations $\pi$ of GSp$(4, \A)$ generated
by $F$ are {\it neither CAP}, short for {\sl Cuspidal Associated to
a Parabolic}, {\it nor endoscopic}, i.e., one arising as the
transfer of a cusp form on GL$(2)\times {\rm GL}(2)$, and
consequently, $L(s,F)$ is not a product of $L$-series of smaller
degrees.

A key example of the first kind of Siegel modular forms we construct
is furnished by the following result, which will be used in [PRa]:

\medskip

\noindent{\bf Theorem A} \, \it Let $E$ be an non-CM elliptic curve
over $\Q$ of conductor $N$. Then there exists a holomorphic, Hecke
eigen-cusp form $F$ of weight $3$ acting on the $3$-dimensional
Siegel upper half space ${\mathfrak H_2}$, such that
$$
L(s, F) \, = \, L(s, {\rm sym}^3(E)).
$$
Moreover, $F$ has the correct level, i.e., what is given by the
conductor $M$ of the Galois representation on the symmetric cube of
the Tate module $T_\ell(E)$ of $E$. \rm

\medskip

Such an $F$ is well known to generate a non-zero $(3,0)$-cohomology
class on a smooth toroidal compactification of the quotient of the
genus $2$ Siegel upper half space ${\mathfrak H}_2$ by a congruence
subgroup $\Gamma$ of level $M$, i.e., which contains the principal
congruence subgroup
$$
\Gamma^{(2)}(M): = \, {\rm Ker}\left({\rm Sp}(4,\Z) \to {\rm
Sp}(4,\Z/M\Z)\right).
$$
When $N$ is {\it square-free}, $M$ turns out to be $N$, and if we
define GSp$(4)$ relative to the alternating form $\begin{pmatrix}0&
I\\-I&0\end{pmatrix}$, we may take $\Gamma$ to be the following:
$$
\Gamma^{(2)}_{I}(M) \, = \, \big\{\gamma = (a_{ij})_{1\leq i,j\leq
4}\in {\rm Sp}(4,\Z) \, \, \vert \, \, a_{ij}\equiv 0\pmod {M}, \,
\forall (i,j)\in \Sigma\big\},
$$
where
$$
\Sigma: = \, \{(12), (31), (32), (41), (42), (43)\}.
$$
Moreover, $F$ is a Siegel newform in the sense of [Sch].
Furthermore, $L(s, {\rm sym}^3(E)) \, = \, \prod_p L_p(s, {\rm
sym}^3(E))$ is the {\it symmetric cube} $L$-function of $E$. To be
precise, if $L(s, E) = \prod_p
\left(1-a_pp^{-s}+p^{1-2s}\right)^{-1}$ is the $L$-series of $E$,
then for any prime $p \nmid N$, $a_p=p+1-\vert E(\F_p)\vert$ equals
$\alpha_p+\beta_p$ with $\alpha_p\beta_p=p$. The Euler factors of
$L(s, {\rm sym}^3(E))$ are given (for $p \nmid N$) by
$$
L_p(s, {\rm sym}^3(E)) \, = \,
\left[(1-\alpha_p^3p^{-s})(1-\alpha^2\beta_pp^{-s})
(1-\alpha_p\beta_p^2p^{-s})(1-\beta_p^3p^{-s})\right]^{-1}.
$$

\medskip

Of course we know now, thanks to the deep works of Wiles, Taylor,
Diamond, Conrad and Breuil ([W], [TW], [BCDT]), that every $E/\Q$ is
modular and attached to a holomorphic newform $f$ on the upper half
plane ${\mathfrak H}_1 \subset \C$ of weight $2$, level $N$ and
trivial character.

Our result extends to give the existence of a holomorphic Siegel
cusp form $F$ corresponding to the symmetric cube of any newform $f$
of even weight $k \geq 2$ and trivial character (see Theorem
$A^\prime$ in Section 1). In particular, newforms for SL$(2,\Z)$ of
weight $2k$, for instance the ubiquitous $\Delta$ function of weight
12, give rise to (vector valued) Siegel modular cusp forms for
Sp$(4,\Z)$. The starting point for the proof of Theorem A is the
{\it symmetric cube transfer}, achieved in [KSh1] by the second
author and H.~Kim, of automorphic forms from GL$(2)$ to GL$(4)$.

\medskip

A typical example of the second kind of holomorphic Siegel modular
forms we construct is given by the following:

\medskip

\noindent{\bf Theorem B} \, \it Let $E$ be an elliptic curve over
$\Q$ and $K$ an imaginary quadratic field with non-trivial
automorphism $\theta$ which does not embed into End$_{\overline
\Q}(E)\otimes \Q$. Pick any algebraic Hecke character $\chi$ of $K$
of even non-zero weight, with associated $\ell$-adic character
$\chi_\ell$, which is anti-cyclotomic, i.e., satisfies $\chi^\theta
= \chi^{-1}$. Then there exists a holomorphic, Siegel modular, Hecke
eigen-cusp form $F$, whose level is the conductor of the
Gal$(\overline \Q/\Q)$-representation on $T_\ell(E)\otimes {\rm
Ind}_K^\Q(\chi_\ell)$, such that, as Euler products over $\Q$,
$$
L(s, F) \, = \, L(s, E_K, \chi).
$$
\rm

\medskip

Here $L(s, E_K, \chi)$ denotes the $\chi$-twisted $L$-function of
$E$ over $K$, which defines a degree $2$ Euler product over $K$.
Viewed as a degree $4$ Euler product over $\Q$, it identifies with
$L(s,T_\ell(E) \otimes {\rm Ind}_K^\Q(\chi_\ell))$. Here ${\rm
Ind}_K^\Q$ is short-hand for the induction of a Gal$(\overline
\Q/K)$-module to one of Gal$(\overline\Q/\Q)$. If $N$ is prime to
the conductor $N^\prime$ of ${\rm Ind}_K^\Q(\chi_\ell)$, which
equals the norm (to $\Q$) of the conductor of $\chi$ times the
absolute value of the discriminant $D$ of $K$, then we have $M =
N^2{N^\prime}^2$.

\medskip

It was A.~Wiles who brought the problem of proving Theorem B to the
attention of the first author a few years ago, and our work on this
article began with his question. It is important to note that there
is no intersection of our result with the elegant theorem of
R.L.~Taylor, M.~Harris and D.~Soudry ([HTS]), associating
holomorphic Siegel modular forms to certain cusp forms on GL$(2)/K$.
Indeed, for their work to apply in our setting, we would need
$\chi^2$ to be fixed by $\theta$, forcing an anti-cyclotomic $\chi$
to have order $4$, hence of weight zero. We expect Theorem B to hold
for $\chi$ of any weight, as long as it is not fixed by $\theta$,
and a proof will presumably come out of a blending, and
strengthening, of the results of Arthur ([A]) and Whitehouse ([Wh]).

\medskip

In Section 2 we establish an extension (see Theorem B$^\prime$) for
newforms $\varphi$ of weight $k\geq 2$ relative to anti-cyclotomic
characters $\chi$ of any imaginary quadratic field $K$ of weight
$w\ne 0$ of the same parity (as $k$). The resulting genus $2$ Siegel
cusp form $F$ turns out to be scalar-valued of (Siegel) weight $k$
exactly when $k=w+2$.

\bigskip

These two kinds of forms can be interpreted, thanks to some recent
results such as [KSh1,2], as special cases of the following general
construction. First we need some preliminaries. An isobaric
automorphic representation $\pi$ of GL$_n(\A)$ is said to be {\bf
algebraic} ([C$\ell$1]) if the restriction of the associated
$n$-dimensional representation
$\sigma(\pi_\infty\otimes\vert\cdot\vert^{(1-n)/2})$ of the real
Weil group $W_\R$ to its subgroup $\C^\ast$ is of the form
$\oplus_{j=1}^n \chi_j$, with each $\chi_j$ algebraic, i.e., of the
form $z\to z^{p_j}\overline{z}^{q_j}$ with $p_j, q_j \in \Z$. One
says that $\pi$ is {\bf regular} iff
$\sigma(\pi_\infty)\vert_{\C^\ast}$ is a direct sum of characters
$\chi_j$, each occurring with {\it multiplicity one}.

Given a prime $\ell$ and a finite set $S$ of places containing
$\ell$ and $\infty$, we will say that a regular, algebraic cusp form
$\pi$ on GL$(n)/\Q$ is of {\it $\ell$-adic Galois type} relative to
$S$ if there is an $n$-dimensional semisimple $\ell$-adic
representation $\rho_\ell$ of Gal$(\overline \Q/\Q)$ such that for a
finite set of places $S$ containing the ramified primes and $\ell$,
$$
L^S(s-w/2, \pi) \, = \, L^S(s,\rho_\ell),
$$
for a suitable ``weight'' $w$. When $\pi$ is selfdual and has a
discrete series component at a prime $r$, say, one knows by Clozel
[C$\ell$2]) that it is of $\ell$-adic Galois type.

\medskip

\noindent{\bf Theorem C} \, \it Let $\pi$ be a cuspidal automorphic
representation of GL$_4(\A)$ of conductor $M$, whose exterior square
$L$-function has a pole, i.e., one has
$$
-{\rm ord}_{s=1}L(s, \pi; \Lambda^2) = 1.
$$
Then the following hold:
\begin{enumerate}
\item[(a)]There exists a
cuspidal automorphic representation $\Pi$ of GSp$_4(\A)$ (of trivial
central character), which is not CAP or endoscopic, such that
$$
L(s,\Pi; r_5)\zeta(s) \, = \, L(s, \pi; \Lambda^2),\leqno(i)
$$
where the first $L$-function on the left is the degree $5$
$L$-function of $\Pi$; moreover, $\Pi$ has a non-zero vector
invariant under a compact open subgroup $K^M=\prod_p K_p^M$ of ${\rm
GSp}(4,\A_f)$, of level $M$, such that for $p$ dividing $M$ at most
once,
$$
K_p^M= K_{p,I}(M) : = \, \big\{k_p \in {\rm GSp}(4,\Z_p) \, \,
\vert \, \, k_p \equiv \begin{pmatrix}\ast & 0 & \ast & \ast\\
\ast & \ast & \ast & \ast\\0 & 0 & \ast & \ast\\ 0 & 0 & 0 &
\ast\end{pmatrix} \pmod {p^{v_p(M)}}\big\}.
$$
\item[(b)]Suppose $\pi$ is in addition regular and algebraic. Then
there is a cusp form $\Pi^h$ on GSp$(4)/\Q$ with $\Pi_\infty$ in the
holomorphic discrete series, with $\Pi^h$ not CAP or endoscopic,
such that the identity of (i) holds with $\Pi^h$ in the place of
$\Pi$.
\item[(c)]If $\pi$ is as in (b) and is moreover of $\ell$-adic Galois
type for some $(\ell, S)$, then
$$
L(s,\Pi^h; r_4) \, = \, L(s,\pi),\leqno(ii)
$$
where the $L$-function on the left is the degree $4$ (spin)
$L$-function of $\Pi^h$.
\end{enumerate}\rm

\medskip

It is useful to note that when $L(s, \pi; \Lambda^2)$ has a pole at
$s=1$, $\pi$ is necessarily selfdual. This is so because
$$
L(s,\pi \times \pi) \, = \, L(s, \pi; \Lambda^2)L(s,\pi; {\rm
sym}^2),\leqno(0.1)
$$
where the symmetric square $L$-function on the right has no zero at
$s=1$ (since $\pi$ is generic, being cuspidal), and the
Rankin-Selberg $L$-function on the left has a pole at $s=1$ iff $\pi
\simeq \pi^\vee$.

In view of the result of Clozel alluded to just before this Theorem,
we see that the stronger identity (ii) (of Theorem C) holds for any
regular, algebraic $\pi$ with a discrete series component at a prime
with its exterior square $L$-function having a pole at $s=1$.

\bigskip

We make use of a number of beautiful results due to others, such as
the works of Taylor, Clozel, Rallis, Cogdell, Piatetskii-Shapiro,
Ginzburg, Harris, Jiang, Soudry, Laumon, and Weissauer, as well as
the joint works of the second author with Kim.

The proofs in this article would have been much shorter had an
unpublished result of Jacquet, Piatetski-Shapiro and Shalika,
concerning the descent of automorphic forms of symplectic type from
GL$(4)$ to GSp$(4)$, been available. To make up for it, we use a
detour and go via SO$(5)$, which makes us lose some information (due
to the difference between the degree $4$ and degree $5$
$L$-functions), recovered later by an argument involving Galois
representations. In the last Section we speculate on how one could
avoid the arithmetic argument if one knew multiplicity one for
Sp$(4)$; it will strengthen Theorem C.

\bigskip

We would like to acknowledge helpful conversations with J.~Cogdell
and D.~Soudry, and thank A.~Wiles for his question regarding Theorem
B, which started this line of investigation. The first author also
thanks J.~Tilouine for the information that the needed papers of
Laumon and Weissauer (on the cohomology of Siegel modular
threefolds) will be appearing in an Ast\'erisque volume in 2006.
Finally, we would like to thank the NSF for supporting this work
through our individual research grants DMS-0402044 (D.R.) and
DMS-0200325 (F.S.).

\vskip 0.2in

\section{Why Theorem C implies Theorem A}

\bigskip

Let $E$ be an elliptic curve over $\Q$ of conductor $N$. Then by
[BCDT], there is a holomorphic newform $h$ of weight $2$, level $N$
and trivial character such that $L(s,E)=L(s,h)$. Let
$A=\R\times\A_f$ be the adele ring of $\Q$, with
$\A_f=\hat{\Z}\otimes\Q$, where $\hat{\Z}:=\lim_n \Z/n\Z \simeq
\prod_p\Z_p$. Then it is well known that such an $h$ generates a
cuspidal automorphic representation $\eta =
\eta_\infty\otimes\eta_f$ of GL$_2(\A) = {\rm GL}(2,\R)\times {\rm
GL}(2,\A_f)$ of trivial central character and conductor $N$ such
that the archimedean component $\eta_\infty$ is in the lowest
discrete series ${\mathcal D}_2$. The associated representation of
the real Weil group $W_\R:=\C^\ast\cup j\C^\ast,$ $j^2=-1,
jzj^{-1}=\overline{z}$ (for all $z\in \C^\ast$), is given, in the
unitary form, by the following:
$$
\sigma_\infty(\eta) \, \simeq \, {\rm
Ind}_{\C^\ast}^{W_\R}\left(\frac{z}{\vert z\vert}\right).\leqno(1.1)
$$
Consequently, its restriction to $\C^\ast$ is given by the direct
sum of characters in the {\it infinity type}:
$$
p_\infty \, = \, \{\frac{z}{\vert z\vert}, \frac{\overline z}{\vert
z\vert}\}.\leqno(1.2)
$$
This is so because $h$ contributes to the $H^1$ of the modular
curve, which has Hodge weight $1$.

Note that $E$ has complex multiplication (CM) by (an order in) an
imaginary quadratic field $K$ iff $\eta$ is automorphically induced
by a Hecke character of $K$, necessarily of weight $1$. Since we
assume that $E$ has no CM, $\eta$ is not dihedral.

Here is a gentle strengthening of Theorem A:

\medskip

\noindent{\bf Theorem ${\bf A^\prime}$} \, \it Let $\eta$ be a
cuspidal, non-dihedral automorphic representation of GL$(2, \A)$
defined by a holomorphic newform $\varphi$ of even weight $k \geq
2$, level $N$ and trivial character. Then there exists a cuspidal
automorphic representation $\Pi$ of GSp$(4,\A)$ of trivial
character, which is unramified at any prime $p \nmid N$, such that
\begin{enumerate}
\item[(a)]$\Pi_\infty$ is in the
holomorphic discrete series, with its parameter being the symmetric
cube of that of $\sigma_\infty$; \, and
\item[(b)]$$ L(s, \Pi) \, = \,
L(s, \eta, {\rm sym}^3).
$$
\end{enumerate}
Moreover, $\Pi$ has a non-zero vector invariant under a compact open
subgroup of GSp$(4,\A_f)$ of level equal to the conductor $M$ of the
symmetric cube of the Galois representation attached to $\varphi$.
When $M$ is square-free, we may take $K$ to be
$$
K_{I}(M) \, = \, \big\{k_f \in {\rm GSp}(4,\hat{\Z}) \, \, \vert \,
\, k_f \equiv \begin{pmatrix}\ast & 0 & \ast & \ast\\ \ast & \ast &
\ast & \ast\\0 & 0 & \ast & \ast\\ 0 & 0 & 0 & \ast\end{pmatrix}
\pmod {M\hat{\Z}}\big\},
$$
\rm

\medskip

If $\eta_p$ is determined at any unramified prime $p$ by an
unordered pair $\{\alpha_p, \beta_p\}$, then the $L$-function on the
right of (b) has an Euler product whose $p$-factors outside $N$ have
the inverse roots:
$$
\{\alpha_p^3, \alpha_p^2\beta_p, \alpha_p\beta_p^2,
\beta_p^3\}.\leqno(1.3)
$$
In fact, at any prime $p$, the local factors such as $L(s, \eta_p;
{\rm sym}^3)$ can be defined as the local factors of the associated
representations $\sigma(\eta_p)$ of the extended Weil group
$W^\prime_{\Q_p}$, which is possible thanks to the local Langlands
correspondence for GL$(n)$ ([HaT], [He]).

\medskip

Now we show why Theorem C implies Theorem $A^\prime$ (and hence
Theorem $A$). By a theorem of the second author with H.~Kim
([KSh1]), asserting that there is an isobaric automorphic
representation sym$^3(\eta)$ of GL$_4(\A)$ such that
$$
L(s, {\rm sym}^3(\eta)) \, = \, L(s, \eta; {\rm sym}^3),\leqno(1.4)
$$
and moreover, the archimedean parameter of sym$^3(\eta)$ is the
symmetric cube of that of $\eta$.

Put
$$
\pi \, = \, {\rm sym}^3(\eta).\leqno(1.5)
$$
Then, as $\eta$ is defined by a holomorphic newform of weight $k$,
its archimedean parameter is, in the unitary form,
$$
\{\left(\frac{z}{\vert z\vert}\right)^{k-1}, \left(\frac{\overline
z}{\vert z\vert}\right)^{k-1}\}.\leqno(1.6)
$$
It follows that the parameter of $\pi_\infty$ is
$$
\{\left(\frac{z}{\vert z\vert}\right)^{3(k-1)}, \left(\frac{z}{\vert
z\vert}\right)^{k-1}, \left(\frac{\overline z}{\vert
z\vert}\right)^{k-1},\left(\frac{\overline z}{\vert
z\vert}\right)^{3(k-1)}\}.\leqno(1.7)
$$
So $\pi$ is algebraic and regular. Moreover, since the central
character $\omega$ of $\eta$ is trivial, it is selfdual (as
$\eta^\vee \simeq \eta\otimes \omega^{-1}$), and so is $\pi = {\rm
sym}^3(\eta)$.

Next, since $\eta$ is generated by a holomorphic newform of weight
$2k\geq 2$ and level $N$, we know by Deligne the existence, for any
prime $\ell$, of an (irreducible) $\ell$-adic representation
$\tau_\ell$ of Gal$(\overline \Q/\Q))$ such that for all primes
$p\nmid N\ell$, the following identity of local factors holds:
$$
L_p(s-(k-1)/2, \eta) \, = \, L_p(s,\tau_\ell).\leqno(1.8)
$$
It follows that sym$^3(\tau_\ell)$ is associated to $\pi$.
Furthermore, one knows (cf. [Ca]) that the conductors of $\eta$ and
$\tau_\ell$ coincide. When this is coupled with the fact that the
symmetric cube transfer $\eta\to{\rm sym}^3(\eta)$ is functorial at
every prime, one deduces the equality of the conductors of
sym$^3(\eta)$ and sym$^3(\tau_\ell)$.

The automorphic representation $\pi$ is cuspidal, as one sees from
Remark 5.9 and Proposition 5.11 of [Sh3]. Here are two other proofs
of the same assertion, both using more the associated $\ell$-adic
representation. The starting point is again that since $\eta$ is
non-dihedral (by hypothesis), [KSh2] implies that $\pi$ is cuspidal
unless the symmetric square sym$^2(\eta)$ of $\eta$ is monomial,
i.e., admits a non-trivial self-twist. But in that case, by the
Tchebotarev density theorem, sym$^2(\tau_\ell)$ would also admit a
non-trivial self-twist, making its restriction to an open subgroup
Gal$(\overline \Q/F)$ reducible, with $F$ cyclic of degree $3$. But
this would contradict a theorem of Ribet ([Ri]), asserting that the
Zariski closure of the image of Gal$(\overline \Q/\Q)$ in GL$(2,
\Q_\ell)$ (under $\tau_\ell$) contains SL$(2, \Q_\ell)$.
Alternatively, we see that the base change of sym$^2(\pi)$ to the
totally real $F$ is not cuspidal, while the associated Galois
representation sym$^2(\tau_\ell)$ is irreducible when restricted to
Gal$(\overline \Q/F)$, contradicting the main result (for $n=3$) of
[Ra3].

To satisfy the hypotheses of Theorem $C$, it remains only to show
that $\pi$ is of {\it symplectic type}, i.e., that the exterior
square $L$-function $L(s, \pi; \Lambda^2)$ has a pole at $s=1$. For
this we appeal to the (well known) identity
$$
L(s,\pi; \Lambda^2) \, = \, L(s, \eta; {\rm
sym}^4))\zeta(s),\leqno(1.9)
$$
which is immediate at the unramified primes $p$. It is known that
the symmetric $4$-th power $L$-function of $\eta$ has no zero at
$s=1$. (One even knows by H.~Kim ([K]) that there is a corresponding
automorphic form sym$^4(\eta)$ on GL$(5)/\Q$, but we do not need to
appeal to it a this point.) Thus the pole of $\zeta(s)$ at $s=1$
induces one of $L(s,\pi; \Lambda^2)$, and Theorem $C$ can be applied
with $\pi={\rm sym}^3(\eta)$. One gets the truth of Theorems $A$ and
$A^\prime$.

\vskip 0.2in

\section{Why Theorem C implies Theorem B}

\bigskip

Let $E$ be an elliptic curve over $\Q$, with associated cusp form
$\eta$ on GL$(2)/\Q$, and let $K$ be an imaginary quadratic field
(with non-trivial automorphism $\theta$), which does not lie in the
$\Q$-endomorphism algebra of $E$. Then $\eta$ is not dihedral
relative to $K$, or in other words, the base change $\eta_K$ (of
$\eta$ to GL$(2)/K$) remains cuspidal. Pick any anti-cyclotomic
Hecke character $\chi$ of $K$ of of even weight $w \ne 0$. We will
use $\chi$ to again denote the associated idele class character, and
use $\chi_\ell$ to denote the corresponding $\ell$-adic character
([Se]) of Gal$(\overline K/K)$. One has
$$
L(s, E_K, \chi) \, = \, L(s, \eta_K \otimes \chi) \, = \, L(s, \eta
\times I_K^\Q(\chi)),\leqno(2.1)
$$
where $I_K^\Q(\chi)$ denotes the (automorphically induced) cusp form
on GL$(2)/\Q$ attached to $\chi$, defined by a holomorphic newform
of weight $w+1$, constructed in this case by Hecke. Since $\chi$ has
a non-zero weight, it cannot be fixed by $\theta$ and so
$I_K^\Q(\chi)$ is indeed cuspidal. The $L$-function on the right of
(2.1) is the Rankin-Selberg $L$-function attached to the pair
$(\eta, I_K^\Q(\chi))$.

\medskip

Here is a gentle strengthening of Theorem B:

\medskip

\noindent{\bf Theorem ${\bf B^\prime}$} \, \it Let $K$ be an
imaginary quadratic field of discriminant $D<0$ and let
$\delta=\delta_{K/\Q}$ denote the quadratic Dirichlet character of
$\Q$ associated to $K$. Consider any cuspidal, automorphic
representation $\eta$ of GL$(2, \A)$ defined by a holomorphic
newform of weight $k \geq 2$, level $N$ and character $\omega$,
which is trivial if $k$ is even and $\delta$ if $k$ is odd. Let
$\chi$ be an anti-cyclotomic Hecke character of weight $w \ne 0$,
such that $w$ and $k$ have same parity. Then there exists a cuspidal
automorphic representation $\Pi$ of GSp$(4,\A)$ of trivial central
character, which is unramified at any prime $p \nmid ND$, such that
\begin{enumerate}
\item[(a)]$\Pi_\infty$ is in the
holomorphic discrete series, with its parameter being the tensor
product of those of $\eta_\infty$ and $I_K^\Q(\chi)_\infty$; \, and
\item[(b)]As Euler products over $\Q$,
$$
L(s, \Pi) \, = \, L(s, \eta_K\otimes \chi).
$$
\item[(c)]$\Pi$ defines a scalar-valued Siegel modular cusp form $F$
when $k=w+2$.
\end{enumerate}
\rm

\medskip

Now we show why Theorem C implies Theorem $B^\prime$ (and hence
Theorem $B$). As $\chi$ is anti-cyclotomic, it can be fixed by
$\theta$ iff it is quadratic, which is impossible as it has non-zero
weight. Hence $\eta_K\otimes \chi$ can be $\tau$-invariant iff it
admits a non-trivial self-twist by $\chi^2$ (=$\chi/\chi^\theta$),
which is again impossible because any such self-twisting character
has to be at most quadratic. Thus $\eta_K\otimes \chi \not\simeq
\left(\eta_K\otimes\chi\right)^\theta$, and so by the base change
theory of Arthur and Clozel ([AC], see also [HH]), it has a cuspidal
automorphic induction
$$
\pi: = \, I_K^\Q(\eta_K\otimes \chi)\leqno(2.2)
$$
to GL$(4)/\Q$. Alternately, one has a cuspidal automorphic
representation $\eta\boxtimes I_K^\Q(\chi)$ of GL$_4(\A)$ (cf.
[Ra1]), such that
$$
L(s, \eta\boxtimes I_K^\Q(\chi) \, = \, L(s,\eta \times
I_K^\Q(\chi))\leqno(2.3)
$$
It follows by the strong multiplicity one theorem that
$$
\pi \, \simeq \, \eta\boxtimes I_K^\Q(\chi).\leqno(2.4)
$$

Unraveling this isomorphism, we see (using (1.6)) that the
archimedean parameter of $\pi$ is given by
$$
\{\left(\frac{z}{\vert z\vert}\right)^{(k-1)+w},
\left(\frac{z}{\vert z\vert}\right)^{(k-1)-w}, \left(\frac{\overline
z}{\vert z\vert}\right)^{(k-1)-w},\left(\frac{\overline z}{\vert
z\vert}\right)^{(k-1)+w}\}.\leqno(2.5)
$$
Since $w$ is non-zero and has the same parity as $k$, $k-1\pm w$
cannot be zero and $k-1+w\ne k-1-w$. In other words, $\pi$ is
regular and algebraic.

We will now show that $\pi$ is of symplectic type, which will imply
in particular that it is selfdual. For this we appeal to the
identity
$$
L(s,\pi; \Lambda^2) \, = \, L(s,{\rm sym}^2(\eta)\otimes
\chi_0\delta)L(s,\omega\otimes
(I_K^\Q(\chi^2))L(s,\omega\chi_0),\leqno(2.6)
$$
where we have used the fact that the central character of
$I_K^\Q(\chi)$ is $\delta$ times the restriction $\chi_0$ of $\chi$
to (the idele classes of) $\Q$, and that
$$
{\rm sym}^2(I_K^\Q(\chi)) \, \simeq \, I_K^\Q(\chi^2) \boxplus
\chi_0.\leqno(2.7)
$$
Now since $\chi$ is anti-cyclotomic,
$$
\chi_0\circ N_{K/\Q} \, = \, \chi\chi^{\theta} \, = \, 1,\leqno(2.8)
$$
implying that $\chi_0$ is $\delta^w$. Since $w$ and $k$ have the
same parity, as does $\omega$, we deduce that $\omega\chi_0$ is
trivial. Hence by (2.6), $L(s,\pi; \Lambda^2)$ is $\zeta(s)$ times a
product of known $L$-functions which do not have a zero at $s=1$.
Thus $\pi$ is symplectic.

To finish proving that Theorem $C$ implies Theorems $B$ and
$B^\prime$, it remains only to check that the holomorphic Siegel
cusp form $F$ which $\pi$ defines is scalar-valued of weight $k$
when $k=w+2$. But the archimedean parameter of a Siegel modular form
of scalar weight $k\geq 2$ is known to be of the form
$$
\{\left(\frac{z}{\vert z\vert}\right)^{2k-3}, \frac{z}{\vert
z\vert}, \frac{\overline z}{\vert z\vert},\left(\frac{\overline
z}{\vert z\vert}\right)^{2k-3}\}.\leqno(2.9)
$$
Now we are done by comparing this with (2.5).

\qed

\vskip 0.2in

\section{Transfer to GSp$(4)$ via SO$(5)$}

\bigskip

Let $\pi$ be as in Theorem $C$. Since by hypothesis, $L(s, \pi;
\Lambda^2)$ has a pole at $s=1$, where $S$ is a finite set of places
containing the ramified and archimedean primes, its global parameter
takes values in Sp$(4,\C)$, which is the $L$-group of SO$(5)$. By
the {\it descent theorem} of Ginzburg, Rallis and Soudry ([GRS]), we
can find a cuspidal, globally generic automorphic representation
${\Pi^{\prime\prime}}$ of the split SO$(5,\A)$ with parameters in
Sp$(4, \C)$.

One knows that SO$(5)$ has a double cover, namely Spin$(5)$, which
is isomorphic to Sp$(4)$. Using this we lift ${\Pi^{\prime\prime}}$
to a cuspidal, globally generic automorphic representation
$\Pi^\prime$ of Sp$(4,\A)$ with trivial central character.

Let $r$ be the standard ($5$-dimensional) representation of the dual
group of Sp$(4)$, which is SO$(5,\C)$.

\medskip

\noindent{\bf Proposition 3.1} \, \it We have
\begin{enumerate}
\item[(a)]
$$
L(s,\Pi^\prime; r)\zeta(s) \, = \, L(s,\pi; \Lambda^2),\leqno(3.1)
$$
\item[(b)]
$$
\sigma_\infty(\Pi^\prime) \oplus 1 \, \simeq \,
\Lambda^2(\sigma_\infty(\pi)).
$$
\end{enumerate}
\rm

\medskip

{\it Proof}. \, By the work of the second author with Cogdell, Kim
and Piatetski-Shapiro ([CoKPSS]), we can transfer $\Pi^\prime$ back
to a cusp form $\pi^\prime$ on GL$(4)/\Q$ such that the arrow
$\Pi_v^\prime \to \pi_v^\prime$ is compatible with the descent of
[GRS], and is functorial at every place $v$ (see [So], [JiSo1]). So
$\pi^\prime$ and $\pi$ are equivalent almost everywhere, hence
isomorphic by the strong multiplicity one theorem. So the
composition of the parameters of ${\Pi^{\prime\prime}}$ with the
natural embedding Sp$(4,\C) \hookrightarrow {\rm GL}(4,\C)$ are the
same as the parameters of $\pi$ at the various places $v$. Next
denote by $\psi$ the pull-back map on the set of automorphic
representations of SO$(5, \A) \simeq {\rm PSp}(4,\A)$, taking values
in the set of automorphic representations of Sp$(4,\A)$, both modulo
equivalence. Then $\psi$ corresponds to the $L$-homomorphism
$$
{}^L\psi: {\rm Sp}(4,\C) \, \rightarrow \, {\rm
SO}(5,\C),\leqno(3.2)
$$
defined so that the exterior square of the standard $4$-dimensional
representation of Sp$(4,\C)$ is isomorphic to $\underline{1}\oplus
{}^L\psi$. The assertions of the Proposition now follow.

\qed

\medskip

The next object is to find a cuspidal, generic automorphic
representation of GSp$(4, \A)$ whose restriction to Sp$(4,\A)$
contains $\Pi^\prime$. This can be done by imitating what Labesse
and Langlands do for SL$(2)$ ([LL]). But we want to refine their
construction in such a way that we keep track of what happens at the
finite primes in order that we do not introduce new ramification.
Here is what we do

First note that in the descent construction $\pi\to \Pi^\prime$ of
Ginzburg, Rallis and Soudry, the authors first define an Eisenstein
series $E(g,s;\pi)$ on SO$(9,\A)$ and take its residue $E_1(g:\pi)$.
Then they define $\Pi^\prime$, which is $\sigma_2(\pi)$ in their
notation, by restricting to SO$(5,\A)$ of a sequence of Fourier
coefficients of $E_1(g;\pi)$ along a unipotent subgroup $U(\A)$ of
SO$(9,\A)$. All this involves a sequence of integrals, and it
follows that if $\pi$ has conductor $M$, then $\Pi^\prime$ admits a
non-zero vector invariant under a principal compact open subgroup
$K^\prime(M)$ of level $M$. To elaborate, here we use the known fact
that the conductor of $\pi$ being $M$ implies that it contains a
non-zero vector fixed by a congruence subgroup of level $M$.

Next we extend $\Pi^\prime$ to a representation $\Pi_1$ of the group
$H:= {\rm Sp}(4,\A)Z(\A)$, where $Z$ denotes the center of GSp$(4)$,
such that $\Pi_1$ is trivial on $Z(\A)$. This makes sense because
$\Pi^\prime$ has trivial central character. Since $\pi$ is
unramified outside (the primes dividing) $N$, $\Pi_1$ is also
unramified outside $N$. Moreover, the transfer at the ramified
primes is also the right one and hence respects the level. Note that
Sp$(4,\A)Z(\A)$ is a normal subgroup of the symplectic similitude
group GSp$(4, \A)$ with a countable quotient group. Now induce
$\Pi_1$ to GSp$(4, \A)$, and choose (as follows) a cuspidal
automorphic representation $\Pi$ of GSp$(4, \A)$ occurring in the
induced representation, which is necessarily globally generic of
trivial central character. Let $K(M)$ denote a principal congruence
subgroup of GSp$(4,\A_f)$ such that $K_1(M):=K(M) \cap {\rm
Sp}(4,\A_f)$ has image $K^\prime(M)$ in SO$(5,\A_f)\simeq {\rm
PSp}(4,\A_f)$. Since $\Pi_1$ has a fixed vector under $K_1(M)$, the
induced representation will, by Frobenius reciprocity, have at least
one constituent $\Pi$ which will have a vector fixed under $K(M)$,
and such a $\Pi$ is what we choose. In particular, $\Pi_p$ is
unramified whenever $\pi$ is.

Suppose a prime $p$ divides $M$. Then it is not hard to see that
$\Pi$ is not unramified, the reason being that the descent of [GRS]
is compatible with the transfer of [CKPSS], which preserves the
epsilon factors. So, in particular, if $p$ divides $M$ exactly once,
then by appealing to [Sch], we see that $\Pi_p$ must in fact contain
a non-zero vector $v$ which is invariant under the {\it Iwahori
congruence subgroup} $K_{p,I}(M)$ defined in the statement of
Theorem C, part (a), in the Introduction, and it is not unramified
for a maximal compact subgroup.

\medskip

For $j=4,5$, let $r_j$ denote the $j$-dimensional representation of
the dual group of GSp$(4)$, which can be identified with GSp$(4,\C)$
itself since it is isomorphic to GSpin$(5)$. (We want to be careful
because some think of $r_4$ as the standard representation, while
the others think that it is $r_5$.) We get, using part (a) of
Prop.3.1,
$$
L(s,\Pi; r_5)\zeta(s) \, = \, L(s,\pi; \Lambda^2).\leqno(3.3)
$$
Indeed, as $\Pi$ has trivial central character, the $r_5$-parameter
of $\Pi_v$ is, at any place $v$, the same as the $r$-parameter of
$\Pi_{1,v}$.

\medskip

\noindent{\bf Proposition 3.4} \, \it Let $\Pi$ be as above, being
associated to $\pi$. Suppose $\pi$ is regular and algebraic. Then
\begin{enumerate}
\item[(a)]$\Pi_\infty$
is of cohomological type;
\item[{}]and
\item[(b)] The restriction to
$\C^\ast$ of the $r_4$-parameter of $\Pi_\infty$ is the same as that
of $\pi_\infty$.
\end{enumerate}
\rm

\medskip

{\it Proof}. \, Since $\pi_\infty$ is regular and algebraic, part
(a) follows from part (b).  Thanks to Proposition 3.1 and the
discussion above, we get
$$
\sigma_\infty(\Pi) \oplus 1 \, \simeq \,
\Lambda^2(\sigma_\infty(\pi)).\leqno(3.5)
$$
Since
$$
\Lambda^2\circ r_4 \, \simeq \, r_5 \oplus \lambda,
$$
where $\lambda$ is the polarization, corresponding to the central
character. Consequently,
$$
\Lambda^2(r_4(\sigma_\infty(\Pi))) \, \simeq \,
\Lambda^2(\sigma_\infty(\pi)).\leqno(3.6)
$$
The parameter of $\pi_\infty$ when restricted to $\C^\ast$ will be
given by a pair $(a,b)$ with $a,b, a+ b \ne 0$ and $a> b$, or more
precisely,
$$
\{\left(\frac{z}{\vert z\vert}\right)^{a}, \left(\frac{z}{\vert
z\vert}\right)^{b}, \left(\frac{\overline z}{\vert
z\vert}\right)^{b},\left(\frac{\overline z}{\vert
z\vert}\right)^{a}\},\leqno(3.7)
$$
The exterior square parameter will then be determined by the pair
$(a+b,a-b)$ (discounting the exponents $0$). Similarly, if the
$r_4$-parameter of $\Pi_\infty$ is given by $(c,d)$ with $c>d$, then
its exterior square will be determined by $(c+d,c-d)$. It follows
that $a=c, b=d$, proving the Proposition.

\qed

\vskip 0.2in

\section{Proof of the first two parts of Theorem C}

\bigskip

Let $\Pi$ be attached to a cusp form $\pi$ of symplectic type, as in
Section 3 above. Then $\Pi$ is a globally generic, cuspidal
automorphic representation of GSp$(4,\A)$ with trivial central
character.

We will henceforth denote by $L(s,\Pi)$ the degree four $L$-function
of $\Pi$.

\medskip

\noindent{\bf Lemma 4.1} \, \it $\Pi$ is not a CAP representation,
nor is it endoscopic (of Yoshida type), i.e., attached to a pair of
cusp forms $\pi_1, \pi_2$ of GL$(2)/\Q$. \rm

\medskip

{\it Proof}. Since $\Pi$ is generic, it cannot be CAP. Suppose it is
endoscopic, associated to a pair of cusp forms $(\pi_1,\pi_2)$ of
GL$(2)/\Q$. Then for any finite set of places containing the
archimedean and ramified ones, we have
$$
L(s,\Pi) \, = \, L(s,\pi_1)L(s,\pi_2).\leqno(4.2)
$$
In this case the central characters of $\pi_1$ and $\pi_2$ are
trivial, and we have by combining (3.1) and (4.2),
$$
L(s,\pi; \Lambda^2) \, = \, L(s,\pi_1 \times
\pi_2)\zeta(s)^2,\leqno(4.3)
$$
where $\pi$ is the cusp form on GL$(4)/\Q$ which descends to $\Pi$.
Since the Rankin-Selberg $L(s, \pi_1\times\pi_2)$ has no zero at
$s=1$, we see that the exterior square $L$-function of $\pi$ has a
double pole. It implies, by (0.1), that $L(s,\pi \times\pi)$ has a
double pole at $s=1$, which violates the cuspidality of $\pi$. Thus
$\Pi$ is not endoscopic.

\qed

\medskip

Part (a) of Theorem C is now proven , and let us now turn our
attention to part (b). From now on, $\pi$ will be regular and
algebraic (in addition to being of symplectic type). Thanks to
Proposition 3.4, and the regularity of $\pi$, we see that
$\Pi_\infty$ is square-integrable. There is a discrete series
$L$-packet $\{\Pi_\infty^h, \Pi_\infty^W\}$ of GSp$(4,\R)$ of
trivial central character, with $\Pi_\infty^h$ being holomorphic and
$\Pi_\infty^W$ being of Whittaker type, i.e., generic, such that
$\Pi_\infty=\Pi_\infty^W$. It is well known that such an
$L$-parameter contributes to the relative Lie algebra cohomology
with coefficients in a finite dimensional representation.
Geometrically, this implies that $\Pi$ contributes to the
$L^2$-cohomology in degree $3$, with coefficients in a suitable
local system $V^{a,b}$, of the complex points of the Shimura variety
$S_K$, defined over $\Q$, attached to a compact open subgroup $K$ of
GSp$(4,\A_f)$. For such a $K$, $\Pi_f^K$ is non-zero and is acted on
by the Hecke algebra ${\mathcal H}({\rm GSp}(4,\A_f), K)$ generated
by the $K$-double cosets; moreover, it determines $\Pi_f$.

By the proof of Zucker's conjecture due to Looijenga and
Saper-Stern, one then obtains, for such a $\Pi$ with $\Pi_f^K\ne 0$,
a class in (the middle intersection cohomology)
$X:=IH^3(S_K^\ast(\C), {\mathcal V}^{a,b})$, where $S_K^\ast$ is the
Baily-Borel-Satake compactification and ${\mathcal V}^{a,b}$ is a
locally constant sheaf attached to $V^{a,b}$. By the comparison
theorem between Betti (singular) and \'etale cohomology, $\Pi_f$
also contributes, for any prime $\ell$, to the $\ell$-adic
cohomology $X_\ell:=IH^3_{\rm et}(S^\ast_K\times\overline\Q,
{\mathcal V}^{a,b}_\ell)$, which is pure (by Gabber) of weight
$a+b-3$, meaning that at any prime $p\ne \ell$ where $X_\ell$ is
unramified as a Galois representation, the inverse roots of the
geometric Frobenius Fr$_p$ are all of absolute value
$p^{(a+b-3)/2}$.

Let $X_\ell(\Pi_f)$ denote the isotypic subspace of
$X_\ell\otimes_{\Q_\ell}\overline\Q_\ell$ defined by $\Pi_f$, which
is again a Gal$(\overline\Q/\Q)$-module since the action of the
Hecke algebra on $X_\ell$ (by algebraic correspondences) commutes
with the Galois action. Its $L$-function is as usual defined by the
Euler product
$$
L(s,X_\ell(\Pi_f)) \, = \, \prod_p \, {\rm det}\big(I-p^{-s}{\rm
Fr}_p\vert X_\ell(\Pi_f)^{I_p}\big)^{-1},
$$
where $I_p$ denotes the inertia group at $p$. The possible Hodge
types of the corresponding summand $X_B(\Pi_f)$ of the singular
cohomology are $(a+b-3,0),(c,d),(d,c),(0,a+b-3)$, with $c+d=a+b-3$.
At the primes $p$ where $X_\ell(\Pi_f)$ is unramified, i.e., where
$I_p$ acts trivially, the Frobenius polynomial (defining the Euler
factor at $p$) has been related to the Hecke polynomial of $\Pi_p$
by Weissauer in his article [Wei], and also by Laumon  in [Lau1,2],
who works instead with the cohomology of $S_K$ with compact
supports. Their beautiful works show, since $\Pi$ is not CAP or
endoscopic, that for some positive integer $m$,
$$
{\rm dim} X_\ell(\Pi_f) \, = \, 4m \, = \,
2(m^W(\Pi_f)+m^h(\Pi_h)),\leqno(4.4)
$$
where $m^W(\Pi_f)$, resp. $m^h(\Pi_f)$, is the multiplicity of
$\Pi_\infty^W\otimes\Pi_f$ ($=\Pi$), resp.
$\Pi_\infty^h\otimes\Pi_f$ in the space of cusp forms on
GSp$(4,\A)$. If either of these representations (of discrete series
type) occurs in the discrete spectrum of $L^2({\rm
GSp}(4,\Q)Z(\A)\backslash {\rm GSp}(4,\A))$, with $Z$ denoting the
center of GSp$(4)$, then it actually occurs in the cuspidal
spectrum. Note that the Hodge type $(a+b-3,0)$ occurs in
$X_B(\Pi_f)$ iff $m^h(\Pi_f)$ is non-zero.

In fact, if $S$ denotes the set of ramified and archimedian primes,
the results of Weissauer and Laumon ({\it loc. cit.}) furnish the
following identity: (with $T=S\cup\{\ell\}$)
$$
L^T(s+(3-a-b)/2,\Pi)^{2(m^W(\Pi_f)+m^h(\Pi_h))} \, = \, L^T(s,
X_\ell(\Pi_f))^4.\leqno(4.5)
$$

\medskip

Since $\Pi$ is globally generic, we know by a theorem of Jiang and
Soudry ([JiSo2]) that $m^W(\Pi_f) =1$. Combining these remarks with
the results of Weissauer and Laumon above, equation (4.4) in
particular, we get
$$
m^h(\Pi_f) > 0.\leqno(4.6)
$$
In other words, the holomorphic representation
$\Pi^h:=\Pi^h_\infty\otimes\Pi_f$ occurs in the space of cusp forms
in $L^2({\rm GSp}(4,\Q)Z(\A)\backslash {\rm GSp}(4,\A))$. Since
$\Pi^h$ is $L$-equivalent to $\Pi$, it has the same $L$-function.
(In fact, the only change is at the archimedean place, where
$\Pi_\infty^h$ and $\Pi_\infty$ are in the same discrete series
$L$-packet.) Hence by appealing to (3.3) and Proposition 3.4, we get
the identity
$$
L(s,\Pi^h; r_5)\zeta(s) \, = \, L(s,\pi; \Lambda^2)\leqno(4.7)
$$
This proves part (b) of Theorem C. Note that (4.5) involves
incomplete $L$-functions, but (4.7) equates complete $L$-functions.

\medskip

\noindent{\bf Remark 4.8} \, {\it Exploiting an idea of Taylor
([Ta]), one can even deduce from Weissauer's theorem that
$m^h(\Pi_f)=1$,} which yields the identity $L^T(s+(3-a-b)/2,\Pi)=
L^T(s, X_\ell(\Pi_f))$. Indeed, thanks to the decomposition theorem
([BBD]), $X_\ell$ can be seen to be a direct factor of $H^3_{\rm
et}({\tilde S}_K\times\overline\Q, {\mathcal V}^{a,b}_\ell)$, where
$\tilde S_K$ denotes a smooth toroidal compactification of $S_K$,
which can be taken to be defined over $\Q$ and such that the
complement of $S_K$ (in $\tilde S_K$) is a divisor $D$ with normal
crossings. It follows, by a theorem of Faltings ([Fa], [Ts]), that
$X_\ell(\Pi_f)$ is of Hodge-Tate type. Let ${\mathcal
G}_\ell(\Pi_f)$ denote the Lie algebra of the Zariski closure of the
image of Galois in GL$(X_\ell(\Pi_f))$. Then by a theorem of Shankar
Sen ([Se]), there is an element $\varphi$ in ${\mathcal
G}_\ell(\Pi_f)$ which has eigenvalues $\alpha, \beta, \gamma,
\delta$, $a+b-3>c>d>0$, with respective multiplicities $m^h(\Pi_f),
m^W(\Pi_f), m^W(\Pi_f), m^h(\Pi_f)$. But then (4.5) forces
$m^h(\Pi_f)= m^W(\Pi_f)$, which is $1$.

\vskip 0.2in

\section{Proof of the last part of Theorem C}

\medskip

We will now prove the part (c) (of Theorem C), under the hypothesis
that $\pi$ is of $\ell$-adic Galois type for some $(\ell, S)$. Let
$\rho_\ell$ be the associated $4$-dimensional, semisimple
$\ell$-adic representation attached to $\pi$, having the same
(degree $4$) $L$-function outside $S$. Let us base change
$\rho_\ell$ to $\overline\Q_\ell$ and work in the absolute setting.
Also, there is a $4$-dimensional Galois representation
$\rho_\ell^\prime$ such that the semisimplification of
$X_\ell(\Pi_f)$ is a multiple of $\rho^\prime_\ell$ (see (4.4),
(4.5)). Then we have two continuous
$\overline\Q_\ell$-representations of Gal$(\overline\Q/\Q)$ which
satisfy, by the Tchebotarev density theorem,
$$
\Lambda^2(\rho_\ell) \, \simeq \,
\Lambda^2(\rho_\ell^\prime).\leqno(5.1)
$$

Note that $\Lambda^2(\rho^\prime_\ell)$ is of orthogonal similitude
type, which can be seen by using the wedge pairing
$$
\Lambda^2(\rho^\prime_\ell) \times \Lambda^2(\rho^\prime_\ell) \,
\rightarrow \, \Lambda^4(\rho^\prime_\ell),\leqno(5.2)
$$
where the right hand side is the determinant of the $4$-dimensional
representation $\rho^\prime_\ell$.

One knows by Urban ([U], Prop. 3.5), that since $m^W(\Pi_f)$ is odd,
the representation $\rho_\ell^\prime$ is of symplectic type, i.e.,
its exterior square representation has an invariant line. In fact,
this can also be deduced by using Remark 4.7, which shows that
$X_\ell(\Pi_f)^{\rm ss} = \rho^\prime_\ell$, and the symplectic
pairing comes from the cup product on the intersection cohomology in
the middle dimension, which is alternating because the dimension is
odd.

Anyhow, since the weight $w$ of $\rho_\ell^\prime$ is $w:=a+b-3$,
its polarization is the $w$-th power of the cyclotomic character
$\chi_{{\rm cyc},\ell}$. We get a decomposition
$$
\Lambda^2(\rho^\prime_\ell) \, \simeq \, r(\rho^\prime_\ell) \,
\oplus \chi_{{\rm cyc},\ell}^{w},\leqno(5.3)
$$
where $r=r_5$ is the $5$-dimensional representation of GSp$(4)$,
taking values in SGO$(5)$, which is the connected component of
GO$(5)$, the orthogonal similitude group; see [Ra2], Section 1, for
a discussion of SGO$(n)$.

Thanks to (5.1), $\rho_\ell$ is also of symplectic type, and we get
a decomposition of $\Lambda^2(\rho_\ell)$, leading to an isomorphism
$$
r(\rho^\prime_\ell) \, \simeq \, r(\rho_\ell).\leqno(5.4)
$$
Now consider the short exact sequence of trivial
Gal$(\overline\Q/\Q)$-modules
$$
1\to\overline\Q_\ell^\ast \to {\rm GSp}(4,\overline\Q_\ell) \to {\rm
SGO}(5,\overline\Q_\ell) \to 1.\leqno(5.5)
$$
The associated long exact cohomology sequence of $\G_\Q={\rm
Gal}(\overline\Q/\Q)$ gives
$$
{\rm Hom}(\G_\Q,\overline\Q_\ell^\ast) \to {\rm Hom}(\G_\Q, {\rm
GSp}(4,\overline\Q_\ell)) \to {\rm Hom}(\G_\Q, {\rm
SGO}(5,\overline\Q_\ell)) \to H^2(\G_\Q,
\overline\Q_\ell^\ast),\leqno(5.6)
$$
where the last group on the right is zero by a theorem of Tate. In
any case, two Galois representations into ${\rm
GSp}(4,\overline\Q_\ell)$ having the same image in ${\rm
SGO}(5,\overline\Q_\ell)$, such as our $\rho_\ell$ and
$\rho^\prime_\ell$, must satisfy
$$
\rho^\prime_\ell \, \simeq \, \rho_\ell\otimes\nu_\ell,\leqno(5.7)
$$
where $\nu_\ell$ is a $1$-dimensional
$\overline\Q_\ell$-representation of $\G_\Q$. Then by (5.1),
$\Lambda^2(\rho_\ell^\prime)$ admits a non-trivial self-twist by
$\nu_\ell^2$. Taking determinants, this implies that
$\nu_\ell^{12}=1$. So $\nu_\ell$ is of finite order and corresponds
to a Dirichlet character $\nu$. Put
$$
\Pi_0^h: = \, \Pi^h\otimes\nu^{-1}.\leqno(5.8)
$$
Now replace $\Pi^h$ by $\Pi_0^h$, we get
$$
L^S(s,\Pi^h; r_4) \, = \, L^S(s,\pi).\leqno(5.9)
$$

It remains to check that (5.9) holds also at the places in $S$. One
knows that for any $v \in S$, the local factor $L(s,\pi_v)$ (resp.
$L(1-s,\pi_v^\vee)$) has no pole in $\Re(s)\geq 1/2$ (resp. in
$\Re(s)\leq 1/2$) (see for instance [BaR]). We {\it claim} that the
same holds for $L(s,\Pi^h_v;r_4)$ (resp.
$L(1-s,{\Pi^h_v}^\vee;r_4)$). Indeed, if $\Pi^h_v$ is tempered, as
it is for $v=\infty$, there is nothing to prove. So we may assume
that $v$ is finite and that $\Pi^h_v=\Pi_v$ is not supercuspidal.
The two maximal parabolic subgroups of GSp$(4)$ are both isomorphic
to GL$(2)\times {\rm GL}(1)$, and by using the knowledge of the
local Langlands correspondence for GL$(2)$ (and GL$(1)$) and a
standard reduction, we see that the analogous correspondence holds
for non-supercuspidal representations of GSp$(4)$. (One also has
this for supercuspidal representations thanks to [JiSo1], but we
will not need it here.) Let $\sigma_v$ be the $4$-dimensional
representation of the extended Weil group $W_{\Q_v}\times {\rm
SL}(2,\C)$ attached to $\Pi_v$. Since $\Pi_v$ is not supercuspidal,
and can moreover be taken to be not square-integrable, $\sigma_v$ is
reducible, and can be written as a sum of two $2$-dimensional
representations $\tau_{1,v}, \tau_{2,v}$ of trivial determinant. On
the other hand, we know that $\Lambda^2(\sigma_v)$ corresponds to
$\Lambda^2(\pi_v)$, and so any one-dimensional representation
occurring in $\Lambda^2(\sigma_v) =\tau_{1,v}\otimes\tau_{2,v}\oplus
1\oplus 1$ will be of norm $< 1/2$, as $\Lambda^2(\pi)$ is (cf. [K])
an isobaric automorphic representation of GL$(6,\A)$. The claim now
follows by considering different possibilities for
$\tau_{1,v},\tau_{2,v}$. Consequently, the poles of
$L(s,\Pi^h_v;r_4)$ (resp. $L(s,\pi_v)$) are disjoint from those of
$L(1-s,{\Pi^h_v}^\vee;r_4)$ (resp. $L(1-s,\pi_v^\vee)$). For any
subset $S_1$ of $S$ consisting of finite primes, choose a character
$\nu=\nu^{S_1}$ which is highly ramified at $S_1$ and trivial at
$S-S_1$. Now we may appeal to [Sh2] to get all the resulting local
factors at the places in $S_1$ to be $1$, without changing the ones
at $S-S_1$. Arguing as in [Ra2], proof of Prop.4.1, we get the
equality $L(s,\pi_v)=L(s,\Pi^h_v;r_4)$ at every $v\in S$.

This finishes the proof of part (c) of Theorem C.

\qed

\vskip 0.2in

\section{Concluding Remarks}

\bigskip

Suppose $\pi, \pi'$ are cuspidal automorphic representations of
GL$(4,\A)$ whose exterior square $L$-functions agree almost
everywhere. It is natural to ask if they must differ by a character
twist. Knowing this {\it a priori} will help remove the hypothesis
in part (c) of Theorem C that $\pi$ is of $\ell$-adic Galois type,
the reason being that we can use the work of the second author with
M.~Asgari ([ASh]) to transfer $\Pi$ back to a cuspidal automorphic
representation $\pi'$ of GL$(4,\A)$, and by construction, $\pi$ and
$\pi'$ will have the same exterior square $L$-function.

When $\pi, \pi'$ are of symplectic type (as in the case of Theorem
C), an affirmative answer should follow from a knowledge of {\it
multiplicity one} for the space of cusp forms on Sp$(4,\A)$. This is
analogous to the situation with SL$(2)$ as in the work of the first
author. To elaborate, one knows multiplicity one for SL$(2)$
([Ra1]), which implies that two cusp forms $\eta, \eta'$ on GL$(2)$
are twist-equivalent if their adjoint $L$-functions are the same
(outside a finite number of factors). Analogously, in our present
case of interest, there should, by functoriality, be cuspidal
automorphic representations $\Pi, \Pi'$ of GSp$(4,\A)$ whose degree
$4$ $L$-functions agree with those of $\pi, \pi'$ respectively. The
hypothesis on the coincidence of the exterior square $L$-functions
of $\pi, \pi'$ implies that $L(s,\Pi; r_5) = L(s,\Pi'; r_5)$. This
should translate, by a study of the stable trace formula of Sp$(4)$,
analogous to that of SL$(2)$ carried out by Labesse and Langlands
([LL]), to the statement that multiplicity one for Sp$(4)$ would
imply the twist equivalence of $\Pi, \Pi'$. Then by the strong
multiplicity one theorem, the same would hold for $\pi, \pi'$.

\medskip

In Theorem B, the weight of the anti-cyclotomic Hecke character
$\chi$ is assumed to be different from $0$. In this (excluded) case,
one can proceed as in this article, but $\Pi_\infty$ will no longer
be in the discrete series, as the parameter of $\pi$ will not be
regular at infinity. It is algebraic and semi-regular, however. It
is likely that $\Pi_\infty$ should be in a (non-degenerate) limits
of discrete series, the reason being tat the polarization of the
associated $4$-dimensional, symplectic representation of
Gal$(\overline\Q/\Q)$ is the transfer of $\chi$ to $\Q$ times the
cyclotomic character, which is {\it odd}. Still, the existence of a
global holomorphic $L$-equivalent $\Pi^h$ is not clear, as one
cannot appeal to the $\ell$-adic cohomology of the Siegel threefold.
However, such a $\Pi^h$ should exist, as predicted by the
conjectures of J.~Arthur, and this should follow from a through
understanding of the comparison of the stable trace formula for
GSp$(4)$ and the twisted, stable trace formula for GL$(4)$ relative
to the automorphism $g\to {}^tg^{-1}$. This appears to be close
given the works of Arthur ([A]) on this comparison, and of
D.~Whitehouse ([Wh]), who recently proved the crucially needed {\it
twisted, weighted fundamental lemma} conjectured by Arthur.

\vskip 0.2in

\section*{\bf Bibliography}

\begin{description}

\item[{[A]}] J.~Arthur, Automorphic representations of ${\rm GSp(4)}$.
in {\it Contributions to automorphic forms, geometry, and number
theory}, 65--81, Johns Hopkins Univ. Press, Baltimore, MD (2004).
\item[{[AC]}] J.~Arthur and L.~Clozel, {\it Simple Algebras, Base
Change and the Advanced Theory of the Trace Formula}, Ann. Math.
Studies {\bf 120} (1989), Princeton, NJ.
\item[{[ASh]}] M.~Asgari
and F.~Shahidi, Generic transfer for general spin groups, Duke Math
Journal {\bf 132}, 137--190 (2006).
\item[{[BaR]}] L.~Barthel and D.~Ramakrishnan,
\emph{A nonvanishing result for twists of $L$-functions of ${\rm
GL}(n)$}, Duke Math. J. {\bf 74},  no. 3, 681--700 (1994).
\item[{[BBD]}] A.~Beilinson, J.~Bernstein and
P.~Deligne, \emph{Faisceaux pervers}, in {\it Analyse et topologie
sur les espaces singuliers I}, Ast\'erisque {\bf 100} (1982),
5--171.
\item[{[BCDT]}]C.~Breuil, B.~Conrad, F.~Diamond and R.L.~Taylor,
\emph{On the modularity of elliptic curves over $\bold Q$: wild
3-adic exercises}, J. Amer. Math. Soc. {\bf 14},  no. 4, 843--939
(2001).
\item[{[BuG]}] D.~Bump and D.~Ginzburg, \emph{Symmetric square $L$-functions on
${\rm GL}(r)$}, Ann. of Math. (2) {\bf 136}, no. 1, 137--205 (1992).
\item[{[Ca]}]H.~Carayol, \emph{Sur les représentations $l$-adiques associées aux formes
modulaires de Hilbert}, Ann. Sci. École Norm. Sup. (4) {\bf 19}, no.
3, 409--468 (1986).
\item[{[C$\ell$1]}] L.~Clozel, \emph{Motifs et formes automorphes},
in {\it Automorphic Forms, Shimura varieties, and $L$-functions},
vol. I, 77--159, Perspectives in Math. {\bf 10} (1990)
\item[{[C$\ell$2]}]L.~Clozel, \emph{Représentations galoisiennes associées
aux représentations automorphes autoduales de ${\rm GL}(n)$}, Publ.
Math. IHES {\bf 73}, 97--145 (1991).
\item[{[CoKPSS]}] J.~Cogdell, H.~Kim, I.~Piatetski-Shapiro and F.~Shahidi, \emph{On lifting
from classical groups to ${\rm GL}\sb N$}, Publ. Math. IHES {\bf 93}
(2001), 5--30.
\item[{[De]}] P.~Deligne, \emph{Les constantes des \'equations fonctionnelles des
fonctions $L$}, in {\it Modular functions of one variable} II,
Springer Lecture Notes {\bf 349} (1973), 501-597.
\item[{[Fa]}] G.~Faltings, \emph{Crystalline cohomology and $p$-adic Galois-representations},
in {\it Algebraic analysis, geometry, and number theory} (1989),
25--80, Johns Hopkins Univ. Press, Baltimore, MD.
\item[{[Fo]}] J.-M.~Fontaine, \emph{Arithm\'etique des
repr\'esentations Galoisiennes $p$-adique}, pr\'epublication d'Orsay
{\bf 24} (March 2000).
\item[{[GJ]}] S.~Gelbart and H.~Jacquet,
\emph{A relation between automorphic
representations of GL$(2)$ and GL$(3)$},
 Ann. Scient. \'Ec. Norm. Sup. (4)
{\bf 11} (1979), 471--542.
\item[{[GRS]}] D.~Ginzburg, S.~Rallis and D.~Soudry, \emph{On a correspondence between
cuspidal representations of ${\rm GL}\sb {2n}$ and $\widetilde{\rm
Sp}\sb {2n}$}, Journal of the AMS {\bf 12} (1999), no. 3, 849--907.
\item[{[HaST]}] M.~Harris, R.L.~Taylor and D.~Soudry,
$\ell$-adic representations associated to modular forms over
imaginary quadratic fields I. Lifting to ${\rm GSp}\sb 4(Q)$,
Invent. Math. {\bf 112},  no. 2, 377--411 (1993).
\item[{[HaT]}] M.~Harris and R.~Taylor, \emph{On the geometry and cohomology of some
simple Shimura varieties}, preprint (2000), to appear in the Annals
of Math. Studies, Princeton.
\item[{[He]}] G.~Henniart, \emph{Une preuve simple des conjectures de Langlands pour
${\rm GL}(n)$ sur un corps $p$-adique}, Invent. Math. {\bf 139}, no.
2, 439--455 (2000).
\item[{[HH]}] G.~Henniart and R.~Herb, \emph{Automorphic induction for ${\rm GL}(n)$
(over local non-Archimedean fields)}, Duke Math. J. {\bf 78}, no. 1,
131--192 (1995).
\item[{[JPSS]}] H.~Jacquet, I.~Piatetski-Shapiro and J.A.~Shalika,
\emph{Rankin-Selberg convolutions}, Amer. J of Math. {\bf 105}
(1983), 367--464.
\item[{[JS1]}] H.~Jacquet and J.A.~Shalika,
\emph{Euler products and the classification of automorphic forms} I
\& II, Amer. J of Math. {\bf 103} (1981), 499--558 \& 777--815.
\item[{[JS3]}] H.~Jacquet and J.A.~Shalika, \emph{Exterior square
$L$-functions}, in {\it Automorphic forms,
Shimura varieties, and $L$-functions}, Vol. II, 143--226,
Perspectives in Math. {\bf 11} (1990), Academic Press, Boston, MA.
\item[{[JiSo1]}]D.~Jiang and D.~Soudry, \emph{Generic representations and local
Langlands reciprocity law for $p$-adic ${\rm SO}\sb {2n+1}$}, in
{\it Contributions to automorphic forms, geometry, and number
theory}, 457--519, Johns Hopkins Univ. Press, Baltimore, MD (2004).
\item[{[JiSo2]}]D.~Jiang and D.~Soudry, \emph{}, \emph{On the multiplicity one
theorem for generic automorphic forms on GSp$(4)$}, preprint (2006).
\item[{[K]}]H.~Kim,  \emph{Functoriality for the exterior square of ${\rm GL}\sb 4$
and the symmetric fourth of ${\rm GL}\sb 2$}, J. Amer. Math. Soc.
{\bf 16}, no. 1, 139--183 (2003).
\item[{[KSh1]}] H.~Kim and F.~Shahidi,
\emph{Functorial products for GL$(2) \times $GL$(3)$ and the
symmetric cube for GL$(2)$}, preprint (2000),  Ann. of Math. (2)
{\bf 155}, no. 3, 837--893 (2002).
\item[{[KSh2]}]H.~Kim and F.~Shahidi, \emph{Cuspidality of symmetric powers with
applications}, Duke Math. J. {\bf 112}, no. 1, 177--197 (2002).
\item[{[KuRS]}] S.~Kudla, S.~Rallis and D.~Soudry, On the degree
$5$ $L$-function for ${\rm Sp}(2)$,  Invent. Math. {\bf 107} (1992),
no. 3, 483--541.
\item[{[LL]}]J.-P.~Labesse and R.P.~Langlands,
\emph{$L$-indistinguishability for ${\rm SL}(2)$}, Canad. J. Math.
{\bf 31}, no. 4, 726--785 (1979).
\item[{[La1]}] R.P.~Langlands, \emph{Problems in the theory of automorphic forms}, in
{\it Lectures in modern analysis and applications III}, Lecture
Notes in Math. {\bf 170} (1970), Springer-Verlag, Berlin, 18--61.
\item[{[La2]}] R.P.~Langlands, \emph{On the notion of an automorphic representation. A supplement},
in {\it Automorphic forms, Representations and $L$-functions}, ed.
by A. Borel and W. Casselman, Proc. symp. Pure Math {\bf 33}, part
1, 203-207, AMS. Providence (1979).
\item[{[Lau1]}] G.~Laumon,  \emph{Sur la cohomologie à supports compacts des vari\'et\'es de Shimura
pour ${\rm GSp}(4)/\Q$}, Compositio Math. {\bf 105} (1997), no. 3,
267--359.
\item[{[Lau2]}] G.~Laumon, \emph{Fonctions z\'etas des vari\'et\'es de Siegel de dimension trois},
in {\sl Formes automorphes II: le cas du groupe GSp(4)}, Edited by
J. Tilouine, H. Carayol,  M. Harris, M.-F. Vigneras, Asterisque {\bf
302}, Soc. Math. France {\it Ast\'erisuqe} (2006).
\item[{[PRa]}] K.~Paranjape and D.~Ramakrishnan, \emph{Modular forms and
Calabi-Yau varieties}, in preparation.
\item[{[Ra1]}] D.~Ramakrishnan, \emph{Modularity of the Rankin-Selberg $L$-series, and
Multiplicity one for SL$(2)$}, Annals of Mathematics {\bf 152}
(2000), 45--111.
\item[{[Ra2]}] D.~Ramakrishnan, \emph{Modularity of solvable Artin
representations of GO$(4)$-type}, International Mathematical
Research Notices (IMRN) {\bf 2002}, No. {\bf 1} (2002), 1--54.
\item[{[Ra3]}] D.~Ramakrishnan, \emph{Irreducibility and cuspidality},
in {\it Representation Theory and Automorphic Forms}, Progress in
Mathematics, Springer-Birkhaeuser, to appear (2006).
\item[{[Ri]}] K.~Ribet, \emph{Galois representations attached to eigenforms with
Nebentypus}, in {\sl Modular functions of one variable, V }, 17--51,
Lecture Notes in Math. {\bf 601}, Springer, Berlin (1977).
\item[{[Sch]}]R.~Schmidt, \emph{On Siegel modular forms of degree $2$ with square-free
level}, preprint (2006).
\item[{[Sen]}] S.~Sen, Lie algebras of Galois groups arising from
Hodge-Tate modules,  Ann. of Math. (2)  {\bf 97}, 160--170 (1973).
\item[{[Se]}] J.-P.~Serre, \it{Abelian $\ell$-adic representations}, With the
collaboration of Willem Kuyk and John Labute, Revised reprint of the
1968 original, Research Notes in Mathematics {\bf 7}, A.K.~Peters
Ltd., Wellesley, MA (1998).
\item[{[Sh1]}] F.~Shahidi, \emph{On the Ramanujan conjecture and the
finiteness of poles for certain $L$-functions}, Ann. of Math. (2)
{\bf 127} (1988), 547--584.
\item[{[Sh2]}] F.~Shahidi, \emph{Twists of a general class of $L$-functions by highly
ramified characters}, Canadian Math. Bulletin {\bf 43}, 380--384
(2000).
\item[{[Sh3]}] F.~Shahidi, \emph{Infinite dimensional groups and
automorphic $L$-functions}, Pure and Applied Math. Quarterly {\bf
1}, No. 3, {\sl Special issue in memory of Armand Borel}, part 2 of
3, 683--699 (2005).
\item[{[So]}]D.~Soudry,  \emph{On Langlands functoriality from classical groups
to ${\rm GL}\sb n$}, in {\it Automorphic forms I}, Ast\'erisque {\bf
298}, 335--390 (2005).
\item[{[Ta]}]R.L.~Taylor, \emph{On the $\ell$-adic cohomology of Siegel threefolds},
Inv. Math. {\bf 114}, pp. 289-310 (1993).
\item[{[TW]}]R.L.~Taylor and A.~Wiles, \emph{Ring-theoretic properties of certain Hecke algebras},
Ann. of Math. (2) {\bf 141},  no. 3, 553--572 (1995).
\item[{[Ts]}]T.~Tsuji, \emph{$p$-adic \'etale cohomology and crystalline
cohomology in the semistable reduction case}, Inv. Math. {\bf 137},
No. 2, 233--411 (1999).
\item[{[U]}]E.~Urban, Selmer groups and the Eisenstein-Klingen ideal,
Duke Math. Journal {\bf 106},  no. 3, 485--525 (2001).
\item[{[Wei]}] R.~Weissauer, \emph{Four dimensional Galois
representations}, in {\sl Formes automorphes II: le cas du groupe
GSp(4)}, Edited by J. Tilouine, H. Carayol,  M. Harris, M.-F.
Vigneras, Asterisque {\bf 302}, Soc. Math. France {\it Ast\'erisuqe}
(2006).
\item[{[Wh]}] D.~Whitehouse, \emph{The twisted weighted fundamental lemma
for the transfer of automorphic forms from GSp$(4)$ to GL$(4)$}, in
{\sl Formes automorphes II: le cas du groupe GSp(4)}, Edited by J.
Tilouine, H. Carayol,  M. Harris, M.-F. Vigneras, Asterisque {\bf
302}, Soc. Math. France {\it Ast\'erisuqe} (2006).
\item[{[W]}]A.~Wiles, \emph{Modular elliptic curves and Fermat's last theorem},
Ann. of Math. (2) {\bf 141},  no. 3, 443--551 (1995).

\bigskip

\end{description}

\vskip 0.3in

\[ \begin{array}{ll}
\mbox{Dinakar Ramakrishnan} & \mbox{Freydoon Shahidi}\\
\mbox{253-37 Caltech} & \mbox{Department of Math., Purdue University} \\
\mbox{Pasadena, CA 91125, USA.} & \mbox{West Lafayette, IN 47907, USA.}\\
\mbox{dinakar@caltech.edu} & \mbox{shahidi@math.purdue.edu}
\end{array}\]

\vskip 0.2in

\end{document}